\documentclass[12pt,reqno,a4paper,twoside]{article}
\usepackage{amsmath,amsthm,amstext,amscd,amssymb,euscript}

\renewcommand{\phi}{\varphi}

\renewcommand{\subset}{\subseteq}
\renewcommand{\emptyset}{\varnothing}

\newcommand{\Zd}{\mathbb Z^d}

\renewcommand{\Pr}{\mathbb P}

\def\1{ {\mathit{1} \!\!\>\!\! I} }

\newcommand{\conn}{\longrightarrow}

\newcommand{\nc}        { \conn  {\hspace{-3.0ex} /} \hspace{1.8ex}   }

\newcommand{\cere}{{\cal CR}}

\newtheorem{theorem}{Theorem}[section]
\newtheorem{definition}[theorem]{Definition}

\newtheorem{lemma}[theorem]{Lemma}

\newtheorem{proposition}[theorem]{Proposition}
\newtheorem{corollary}[theorem]{Corollary}

\newcommand{\beq}{\begin{eqnarray}}
\newcommand{\eeq}{\end{eqnarray}}

\newcommand{\be}{\begin{equation}}
\newcommand{\ee}{\end{equation}}

\newcommand{\eq}{\begin{equation}}
\newcommand{\en}{\end{equation}}

\newcommand{\bl}{\begin{lemma}}
\newcommand{\el}{\end{lemma}}

\newcommand{\br}{\begin{remark}}
\newcommand{\er}{\end{remark}}

\newcommand{\bt}{\begin{theorem}}
\newcommand{\et}{\end{theorem}}

\newcommand{\bd}{\begin{definition}}
\newcommand{\ed}{\end{definition}}

\newcommand{\bp}{\begin{proposition}}
\newcommand{\ep}{\end{proposition}}

\newcommand{\bc}{\begin{corollary}}
\newcommand{\ec}{\end{corollary}}

\newcommand{\bpr}{\begin{proof}}
\newcommand{\epr}{\end{proof}}

\newcommand{\bi}{\begin{itemize}}
\newcommand{\ei}{\end{itemize}}

\newcommand{\ben}{\begin{enumerate}}
\newcommand{\een}{\end{enumerate}}

\newcommand{\shift}{\!\!\!\!}

\newcommand{\Z}{\mathbb Z}
\newcommand{\R}{\mathbb R}
\newcommand{\N}{\mathbb N}

\newcommand{\E}{\mathbb E}

\newcommand{\pee}{\ensuremath{\mathbb{P}}}

\newcommand{\ce}{\ensuremath{\mathcal{C}_{\rm le}}}
\newcommand{\cee}{\ensuremath{{\bf \mathcal{C}}}}

\newcommand{\fe}{\ensuremath{\mathcal{F}}}
\newcommand{\mee}{\ensuremath{\mathcal{M}}}

\newcommand{\aaa}{\ensuremath{\mathcal{A}}}

\newcommand{\La}{\ensuremath{\Lambda}}
\newcommand{\la}{\ensuremath{\Lambda}}
\newcommand{\si}{\ensuremath{\sigma}}
\newcommand{\om}{\ensuremath{\omega}}
\newcommand{\epsi}{\ensuremath{\epsilon}}

\newcommand{\QED}{\hspace*{\fill}$\Box$\medskip}

\def\now{
\ifnum\time<60
          12:\ifnum\time<10 0\fi\number\time am
          \else
            \ifnum\time>719\chardef\a=`p\else\chardef\a=`a\fi
          \hour=\time
          \minute=\time
          \divide\hour by 60 
          \ifnum\hour>12\advance\hour by -12\advance\minute by-720 \fi
          \number\hour:%
          \multiply\hour by 60 
          \advance\minute by -\hour
          \ifnum\minute<10 0\fi\number\minute\a m\fi}
\newcount\hour
\newcount\minute
\numberwithin{equation}{section}         

\theoremstyle{remark}

\def\t{{\bf t}}  
\begin{document}
\title{{\bf Maximal clusters in non-critical percolation and
related models}}

\author{
Remco van der Hofstad\footnote{Department of Mathematics and Computer Science,
Eindhoven University of Technology and EURANDOM, P.O.\ Box
513, 5600 MB Eindhoven, The Netherlands.
E-mail: {\tt rhofstad@win.tue.nl, fredig@win.tue.nl}}
\and
Frank Redig$~^*$
}

\maketitle

\footnotesize
\begin{quote}
{\bf Abstract:} We investigate the maximal non-critical cluster
in a big box in various percolation-type models.
We investigate its typical size, and the fluctuations
around this typical size. The limit law of these
fluctuations are related to maxima of {\it independent}
random variable with law described by a single cluster.

\end{quote}
\normalsize

\vspace{12pt}

\section{Introduction and main results}
\label{sec-intro}
Bazant in \cite{baz} studies the distribution
of maximal subcritical clusters, both numerically
and via a non-rigorous renormalization
group argument. He finds that the cardinality of
maximal clusters behaves like the maximum of
independent geometrically distributed
random variables, i.e., a ``Gumbel-like"
distribution. In his paper, the role
of the FKG inequality, which means that clusters ``repel each other",
is already emphasized in a subadditivity argument.

In this paper, we rigourously prove these claims
for a broad class of non-critical percolation type models.
In the FKG context, we can deal both with maximal
subcritical and finite supercritical clusters, and
obtain a Gumbel distribution for both. In a more
general context, we can deal with dependent
percolation models dominated by subcritical
Bernoulli percolation.

The key ingredient of the proof of the Gumbel law is
to use the exponential law for the occurrence
time of rare patterns. This idea is
used by Wyner in \cite{Wyner} in the context of
matching two random sequences. If a cluster
bigger than $u_n$ appears in a box $[-n,n]^d\cap\Zd$ of
volume $(2n+1)^d$, then evidently the occurrence
time $\t_{u_n}$ of such a cluster is less
than $(2n+1)^d$. Therefore,  if $\t_{u_n}$ has
approximately an exponential distribution,
then the probability of having
a cluster larger than $u_n$ is
approximately $1- e^{-(2n+1)^d\pee (\cee_{u_n})}$, where
$\cee_{u_n}$ denotes the event that the
cluster of the origin has cardinality at least
$u_n$. If one can find a scale $u_n = u_n(x)$ such that
$\pee (\cee_{u_n(x)})\simeq e^{-x}/(2n+1)^d$, then one
obtains the Gumbel law. Assuming an exponential
decay of the cluster cardinality, as expected for subcritical percolation, one obtains
$u_n(x) = u_n+ x$, where $u_n=c\log n + o(\log n)$. For finite supercritical
clusters, under the assumption of Weibull-tails where the tails decay
as a stretch exponential with exponent $\delta<1$, we have
$u_n (x) = (c\log n + c'\log\log n + x )^{1/\delta}$.

\subsection{The model}
We consider site percolation and related models
on the lattice $\Z^d$.
A configuration of occupied and vacant sites
is an element $\omega\in\Omega= \{ 0, 1\}^{\Z^d}$.
A site $x$ with $\omega (x)=1$ is called occupied, and a
site with $\omega(x) =0$ is called vacant.


The configuration $\omega$ will be distributed according
to a translation invariant probability measure $\pee$ on the
Borel-$\si$-field of $\Omega$. Examples of $\pee$
include the Bernoulli product measure $\pee_p$ with
$\pee_p(\omega(x) =1) =p$, but we will also consider dependent
random fields, such as the Ising model, below.

A set $A\subset\Z^d$ is {\it connected} if for
any $x,y\in A$ there is a nearest-neighbor path
$\gamma$ joining $x$ and $y$.
The cluster $\cee (x)=\cee (x,\omega)$ of an occupied site $x$
is the largest connected subset of occupied sites
to which $x$ belongs. By convention, $\cee (x) =\emptyset$
if $\omega(x)=0$. We will also need the cluster $\ce (x)$ defined
as follows
    \be
    \ce(x) =
    \begin{cases}
    \cee (x)& \ \text{if}\ x\ \text{is the left endpoint of}\ \cee(x),\\
    \emptyset&\ \text{otherwise.}
    \end{cases}
    \ee
Here by the left-endpoint of a finite set $A\subset\Z^d$, we mean the minimum of $A$ in the
lexicographic order.
By definition $\ce(x)\cap \ce(y) =\emptyset$ if $x\not= y$. In this paper, we will
work with site percolation. In the percolation community,
it is more usual to consider bond percolation (see e.g.\
\cite{Grim99}). However, site percolation is more general than bond percolation,
as shown e.g.\ in \cite[Section 1.6]{Grim99}. We will use results from \cite{Grim99}
proved for bond percolation, but in general these results also hold for site percolation
(as noted in \cite[Section 12.1]{Grim99}).

Percolation has a {\it phase transition}, i.e., for $d\geq 2$,
there exists a critical value $p_c\in (0,1)$ such that
there exists an infinite cluster a.s.\ for $p>p_c$, while no
such cluster exists a.s.\
for $p<p_c$. The goal of this paper is to investigate maximal
clusters in a finite box for $p\neq p_c$.


\subsection{Main results for site percolation}
\label{sec-res}
In this section, we describe our results in the simplest case, namely
for site percolation, where all vertices are independently occupied with probability $p$
and vacant with probability $1-p$.

We will study the maximal cluster inside a big box. To be able to state our result,
we need some further notation. Let $B_n=[-n,n]^d\cap\Z^d$ be the cube of width $2n+1$.
We let
    \be
    \omega_{B_n}
    =
    \begin{cases}
    \omega (x)& \ \text{if}\ x\in B_n,\\
    0&\ \text{otherwise},
    \end{cases}
    \ee
and
    \be
    \label{meenzb}
    \mee_n=\mee_n (\omega) =\max_x|\ce (x,\omega_{B_n})|.
    \ee
The random variable $\mee_n$ is the maximal cluster
inside $B_n$, with zero boundary conditions, i.e., where we do not consider connections outside $B_n$.
The goal of this paper is to obtain an extreme value theorem such as
    \be\label{exval}
    \pee \left(\mee_n \leq u_n +x\right) = e^{-a_n e^{-x}}+o(1)
    \ee
for some $u_n\uparrow\infty$, and
where $a_n$ is a bounded sequence. In words, this means
that the distribution of the maximal cluster is
``Gumbel-like", i.e., looks like the maximum of
independent geometric random variables. The presence
of the bounded sequence $a_n$ is typical for the law of
the maximum of independent geometric random variables, where we do not have
an exact limiting extreme value distribution cannot (see
e.g.\ \cite[Corollary 2.4.1]{G}).

The idea developed in  this paper is that
for any non-critical $p$, the law of $\mee_n$ is asymptotically
equal to the law of the maximum of $(2n+1)^d$ {\it independent}
copies of a random variable $X$ with law
    \eq
    \label{Xdef1}
    \pee(X=n) = \frac 1n \pee (|\cee(0)|=n),
    \en
for $n\geq 1$, and
    \eq
    \label{Xdef2}
    \pee(X=0)= 1-\E(|\cee(0)|^{-1}).
    \en
The law of $X$ in (\ref{Xdef1}--\ref{Xdef2}) turns out to be
equal to the law of the random variable $|\ce(0)|$
(see Lemma \ref{clust} below). Therefore,
the law of $\mee_n$ is equal to the law of $\max_x|\cee (x,\omega_{B_n})|$,
and thus the philosophy of the paper is to show that the clusters are only
weakly dependent. We further use properties of the law of $\cee(0)$ to derive
the asymptotics of $\mee_n$ in more detail.



We note that the cluster size distribution will play an essential
part throughout the proof. We will now state the results on
this cluster size distribution which we need, in order to specialize
the results. Since this law is crucially different for
$p<p_c$ and $p>p_c$, we distinguish these two cases.

For $p<p_c$,  it is shown in \cite[Theorem (6.78)]{Grim99} that
    \eq
    \label{zetadef}
    \zeta(p,d) =  \lim_{n\rightarrow \infty} -\frac 1n \log \pee_p(|\cee(0)|\geq n)
    \en
exists, and that $\zeta=\zeta(p,d)>0$ for all $p<p_c$. Moreover,
there exists $C=C(p)$ such that
    \eq
    \label{CSbdsub}
    \pee_p(|\cee(0)|= n)\leq Cn e^{-\zeta n}.
    \en

We will sometimes work under an assumption that a somewhat stronger
version of (\ref{zetadef}) holds, namely that
    \eq
    \label{subass}
    \lim_{n\rightarrow \infty} \frac{\pee_p(|\cee(0)|\geq n+1)}{\pee_p(|\cee(0)|\geq n)}
    =e^{-\zeta}.
    \en
Assumption (\ref{subass}) is stronger than (\ref{zetadef}), but strictly
weaker than the widely believed tail-behavior, namely that there exist $\theta=\theta(d)\in \R$ and $A=A(p,d)$
such that
    \eq
    \pee_p(|\cee(0)|\geq n) = A n^{\theta} e^{-\zeta n} [1+o(1)].
    \label{conj}
    \en

Our main result for $p<p_c$ is the following theorem:
    \bt
    \label{thm-Gumbel}
    Fix $p<p_c$ and assume that (\ref{subass}) holds.
    Then there exists a sequence $u_n\in \N$, with $u_n\to\infty$, a real number $a>0$ and
    a bounded sequence $a_n\in [a,1]$, such that for all $x\in\N$
        \be
        \pee (\mee_n \leq u_n + x)=e^{-a_n e^{-x\zeta}}+o(1).
        \ee
    \et
Theorem \ref{thm-Gumbel} shows that $\mee_n$ is bounded above and
below by Gumbel laws, and shows in particular that the sequence
$\mee_n-u_n$ is tight. Our proof will reveal that
Theorem \ref{thm-Gumbel} can be extended to yield weak
convergence along certain exponentially growing sequences.

We now go to supercritical results. Since $p_c(1)=1$, we may assume that we are
in dimension $d>1$. When $p>p_c$, then it is shown
in \cite[Theorem (8.61) and (8.65)]{Grim99} that there exist $\eta=\eta(p,d)$
and $\gamma=\gamma(p,d)$ such that
    \eq
    \label{CSbdsup}
    e^{-\gamma n^{\frac{d-1}{d}}} \leq \pee_p(|\cee(0)|= n)\leq e^{-\eta n^{\frac{d-1}{d}}}.
    \en
In $d=2,3$, it is known that the limit
    \eq
    \label{etadef}
    \eta(p,d) =  \lim_{n\rightarrow \infty} -\frac 1{n^{\frac{d-1}{d}}}
    \log \pee_p(n\leq |\cee(0)|<\infty).
    \en
exists. The limit in (\ref{etadef}) is related to the large deviations of large
finite supercritical clusters, and can be written explicitly as a variational problem
over possible cluster shapes. This variational problem involves the surface tension, and
is maximized by the so-called Wulff shape. The result in $d=2$ is in \cite{ACC90, Cerf99},
while for $d=3$, it is in \cite{Cerf00}.

We will again formulate a stronger version of (\ref{etadef}), namely that
for every $x\in \R$, we have
    \eq
    \label{supass}
    \lim_{n\rightarrow \infty} \frac{\pee_p(n+xn^{1/d}\leq |\cee(0)|<\infty)}
    {\pee_p(n\leq |\cee(0)|<\infty)}
    =e^{-x\eta\frac{d-1}{d}},
    \en
and change the definition of $\mee_n$ slightly to
    \be
    \label{meenzbsup}
    \mee_n=\mee_n (\omega) =\max_{x: |\ce(x)|<\infty} |\ce
    (x,\omega_{B_n})|,
    \ee
i.e., we take the largest {\it finite} cluster.
Of course, for $p<p_c$ (\ref{meenzb})
and (\ref{meenzbsup}) coincide.

 Then we can prove
the following scaling property:

    \bt
    \label{thm-Gumbelsup}
    Fix $p>p_c$ and assume that (\ref{supass}) holds.
    Then there exists a sequence $u_n\to\infty$ with
    such that for all $x\in \R$
        \be
        \lim_{n\rightarrow \infty} \pee (\mee_n \leq u_n + x u_n^{1/d})=e^{-e^{-x\eta\frac{d-1}{d}}}.
        \ee
    \et

Theorems \ref{thm-Gumbel} and \ref{thm-Gumbelsup} study fluctuations
of $\mee_n$ around their asymptotic mean under the Assumptions (\ref{subass})
and (\ref{supass}). The main difference between Theorems \ref{thm-Gumbel}
and \ref{thm-Gumbelsup} is that Theorem \ref{thm-Gumbelsup} implies
weak convergence of the rescaled $\mee_n$ since the fluctuations
grow with $n$, whereas in Theorem \ref{thm-Gumbel} this weak
convergence does not hold due to the fact that the fluctuations are
of order 1 so that the discrete nature of cluster sizes persists.

In Section 3 below, we will formulate more
general results that hold {\it without} Assumptions (\ref{subass})
and (\ref{supass}), but that take a form which is less elegant.
It is not so hard to see that one can choose
    \eq
    u_n = O(\log{n})
    \en
for $p<p_c$, while
    \eq
    u_n = O((\log{n})^{\frac{d}{d-1}})
    \en
when $p>p_c$. From Theorems \ref{thm-Gumbel} and \ref{thm-Gumbelsup}
it immediately follows that $\mee_n$ divided by $\log{n}$ for $p<p_c$, respectively,
$(\log{n})^{\frac{d}{d-1}}$  for  $p>p_c$, converges
in probability to a constant. In the next theorems, we will investigate
the typical size of $\mee_n$ in more detail and prove convergence almost surely.

\bt
\label{cor-sub} For $p<p_c$,
    \eq
    \frac{\mee_n}{\log{n}} \rightarrow d\zeta(p,d)\qquad a.s.
    \en
\et

\bt
\label{cor-sup} For $p>p_c$, and $d=2,3$,
    \eq
    \label{maxscal}
    \frac{\mee_n}{(\log{n})^{\frac{d}{d-1}}}
    \rightarrow d^{\frac{d-1}{d}} \eta(p,d)\qquad a.s.
    \en
For $d\geq 4$, if the limit in (\ref{etadef}) exists, then (\ref{maxscal}) holds.
\et

We close this section with a few observations concerning the role of the
boundary conditions. In (\ref{meenzb}), we have taken the maximal cluster
under the zero boundary condition, so that we can write
$\mee_n=\mee_n^{\scriptscriptstyle {\rm (zb)}}$. Alternatively, we could defined
$\mee_n$ under free boundary conditions, i.e.,
    \be
    \label{meenfb}
    \mee_n^{\scriptscriptstyle {\rm (fb)}} =\max_{x\in B_n}|\cee (x,\omega)|,
    \ee
or under periodic boundary conditions, i.e.,
    \be
    \label{meenpb}
    \mee_n^{\scriptscriptstyle {\rm (pb)}} =\max_{x\in B_n}|\cee (x,\omega_{B_n}')|,
    \ee
where $\omega_{B_n}'$ is the site percolation configuration on the torus with vertex set
$B_n$. We will finally show that this makes no difference whatsoever:

    \bt
    \label{thm-bc}
    For $p<p_c$,
        \be
        \pee_p (\mee_n^{\scriptscriptstyle {\rm (zb)}}\neq \mee_n^{\scriptscriptstyle {\rm (fb)}})=
        o(1),
        \qquad \pee_p (\mee_n^{\scriptscriptstyle {\rm (zb)}}\neq \mee_n^{\scriptscriptstyle {\rm (pb)}})=o(1).
        \ee
    For $p>p_c$,
        \be
        \pee_p (\mee_n^{\scriptscriptstyle {\rm (zb)}}\neq \mee_n^{\scriptscriptstyle {\rm (fb)}})=o(1).
        \ee
    \et
Theorem \ref{thm-bc} immediately shows that all results proved for
$\mee_n^{\scriptscriptstyle {\rm (zb)}}$ immediately also apply to
$\mee_n^{\scriptscriptstyle {\rm (fb)}}$ and $\mee_n^{\scriptscriptstyle {\rm (pb)}}$
for $p<p_c$ and to $\mee_n^{\scriptscriptstyle {\rm (fb)}}$ for
$p>p_c$i.e., that the boundary condition is irrelevant.
For $p>p_c$, $\mee_n^{\scriptscriptstyle {\rm (pb)}}$ is
more difficult to work with since it is harder to `exclude' the
infinite cluster on the torus without looking outside the torus.

\subsection{Discussion of the results}
In this section, we discuss our results and their
relation to the literature.

\subsubsection{Runs and one-dimensional site percolation}
In the case where $d=1$, it easily follows that for any $p<p_c=1,$
    \eq
    \pee(|\ce(0)| \geq n) = p^n.
    \en
In this simple case, the largest cluster is equal to the longest
run of ones in $n$ independent tosses. This is a classical problem, and
the leading order asymptotics $\mee_n=\log{n}/\log{p}(1+o(1))$ is the celebrated
Erd\"os-R\'enyi law \cite{ER70}. Our results studies fluctuations around the
Erd\"os-R\'enyi law. This problem has attracted considerable
attention due to its relation to matching problems arising in
sequence alignment (see e.g., \cite{KD92} and the references therein).

\subsubsection{Results for general subcritical FKG models and related Gumbel laws}
Our results for subcritical clusters hold more generally than just for independent site percolation.
The main technical ingredient in the proof are the FKG-inequality and
bounds on the tails of the cluster size distribution. In Section 3 below,
we will state a general result, that can be
proved for site percolation and applied in the context of the
following examples.

\ben
\item The two-dimensional Ising model at $\beta <\beta_c$.

\item The Ising model in general
dimension, at high temperature and/or
high enough magnetic field (see \cite{ghm}).

\item Gibbs measures where the potential has a sufficiently small
Dobrushin norm and a sufficiently high magnetic field.

\een

See \cite{ghm} for an introduction of the Ising model and Section 3 for more details.

We expect that related results hold for other maximal values of
cluster characteristics. Examples are the maximal diameter of a supercritical
finite cluster, or the maximal occupied line (i.e., a sequence of bonds)
with any orientation for $p<1$.
We also expect that our results for maximal finite supercritical clusters continue
to hold in the context of
the Ising model in dimensions $d=2,3$ for $\beta>\beta_c$,
where the Wulff crystal has been identified (see e.g. \cite{BIV00, CP01}), and hence the
exact behavior of the cluster tail is known. We use the version of the exponential
law from \cite{acrv}, and this does not hold for the low-temperature Ising model.
The weaker version of the exponential law proved in \cite{cr} does apply, but
it is unclear whether we can apply this result in the present setting. The proofs
of Theorems \ref{cor-sub}--\ref{cor-sup} are more robust, and are likely to apply
to the Ising model as well.


\subsubsection{Maximal clusters for critical percolation}
Our results are only valid for non-critical percolation. In critical
percolation, the behavior of the largest cluster in a box should be
entirely different. Firstly, the scaling of the largest cluster in
a box should be {\it polynomial} in the volume of the box, rather
then {\it polylogarithmic} as in Corollaries \ref{cor-sub} and \ref{cor-sup}.
Secondly, when properly rescaled, the size of the largest cluster should
converge to a proper random variable, rather than to a constant as in
Corollaries \ref{cor-sub} and \ref{cor-sup}. Thirdly, we expect that
in some cases, the size of the largest cluster depends on the boundary
conditions, which is not true off the critical point (see Theorem
\ref{thm-bc}).

There have been results in the direction of the above claims.
In \cite{BCKS01}, the largest critical cluster in a box was investigated and,
under certain scaling assumptions, it was proved that the largest cluster
with zero boundary conditions scales like $n^{\frac{d\delta}{\delta+1}}$,
where $\delta$ is the critical exponent related to the critical cluster
distribution
    \eq
    \frac{1}{\delta} = \lim_{n\rightarrow \infty}
    -\frac{\log{\pee_{p_c}(|\cee(0)|\geq n)}}{\log n}.
    \en
Of course, it is not obvious that this limit exists. The scaling assumptions
are not expected to be true above the critical dimension $d_c=6$. In \cite{BCHSS04a},
it is conjectured that the same should be true for the largest cluster with
periodic boundary conditions, i.e., for percolation on the torus.

Above the critical dimension $d_c=6$, other scaling occurs. Aizenman \cite{Aize97}
proves that, under a certain assumption on the two-point function, that
the largest cluster has size $n^4$, and that
there are $n^{d-6}$ clusters of such order. The assumption was proved to hold for
nearest-neighbor bond percolation in sufficiently high dimension in \cite{Hara00},
and for sufficiently spread-out percolation above 6 dimensions in \cite{HHS03}.
For periodic boundary conditions, Aizenman \cite{Aize97} conjectured that the scaling
should be like $n^{2d/3}$.
Partial results in this direction have appeared in \cite{BCHSS04a, BCHSS04b}.
It is a well-know result that for the critical random graph, where $p=1/N$,
and $N$ the size of the graph, the largest cluster is of the order $N^{2/3}$.
Thus, Aizenman's conjecture amounts to the conjecture that the largest critical
cluster with periodic boundary conditions scales like the largest critical
cluster for the random graph (see also Section \ref{sss:RG} below).

\subsubsection{Relation to random graphs}
\label{sss:RG}
There is a wealth of related work for random graphs, which are finite
graphs where edges are removed independently. This research topic was
started by a seminal paper of Erd\"os and R\'enyi \cite{ER61}, which
created the field of random graphs. Erd\"os and R\'enyi investigate
what is called {\it the} random graph, i.e., the complete graph
where edges are kept independently with fixed probability $p$ and
removed otherwise. See the books \cite{AS00, Boll01, JLR00}
and the references therein. The fields of percolation and random graphs
have to a large extent evolved independently, with different terminology
and methodology. Only recently, attempts have been made to use the developed
methodology in the other fields (see e.g.\ \cite{BCKS99, BCKS01, BCHSS04a, BCHSS04b,
BCHSS04c}). When dealing with random graphs, it is natural to investigate the
largest connected component or cluster when the size of the graph tends to\
infinity. Therefore, results such as the ones presented in Section \ref{sec-res}
have appeared in this field. In particular, detailed estimates
of large subcritical clusters and supercritical cluster have been obtained.
Of course, for {\it finite} graphs, it is already non-trivial to define what
a critical value is. Above the critical value, the largest cluster has
a size of order of the size of the graph, while below the critical value, the
largest cluster is logarithmic in the size of the graph.

In random graph theory, often there is a {\it discrete duality principle}, which means
that when we remove the largest supercritical cluster, then the size
and distribution of the remaining clusters is very much alike the
size and distribution of subcritical clusters. See e.g.\ \cite[Section 10.4]{AS00}
for an explanation of this principle for branching processes as well
as for the random graph. We note that this principle is {\it false} for
site percolation, as Theorems \ref{cor-sub} and \ref{cor-sup} show.
This distinction arises from the fact that the random graph has no
{\it geometry}, whereas the geometry is essential in the description of large
finite supercritical clusters and appears prominently in the Wulff
shape. It would be of interest to apply our methods
to random graphs.

\subsubsection{Organization}
Our paper is organized as follows. In Section
2, we give heuristics for our results.
In Section 3, we state our general results
for FKG models under certain conditions.
Section 4 is devoted to the proofs of the main
results.

%
%
%
%
%

\section{Extremes and rare events: heuristics}
\label{sec-heur}
We will be interested in the cardinality of maximal clusters
inside a big box. Recall that $B_n=[-n,n]^d\cap \Z^d$.
For $n\in\N$, define the $\si$-field $\fe_n= \fe_{B_n}$.
A {\em pattern} $A_n$ is a configuration with support
on $B_n$, i.e., it is  an element of $\{ 0, 1\}^{B_n}$. We
will identify a pattern with its cylinder, i.e., we
will also denote $A_n$ to be the set of those $\omega$ such
that $\omega_{B_n} = A_n$. For a pattern $A_n$, we define
its occurrence time to be
    \be
    \t_{A_n} (\omega) = \min \Big\{ |B_k|: \exists x\in\ B_k
    \ {\rm such~that~}\ B_n+x\subset B_k \ {\rm and} \ \theta_{x}\omega_{B_n},
    =  A_n \Big\}
    \ee
where $\theta_x \omega$ denotes the configuration $\omega$ shifted over $x$, so that
$(\theta_x \omega)(y)=\omega(x+y)$.
In words, this is the volume of the minimal cube $B_k$
which ``contains" the pattern $A_n$. One expects
that $\t_{A_n}$ is of the order $\pee(A_n)^{-1}$.
For $E_n\in\fe_n$, there exists a unique
set of patterns $\aaa (E_n)$ such that:
    \[
    E_n = \bigcup_{A_n\in\aaa (E_n)} A_n
    \]
The occurrence time of $E_n$ is then defined
as
    \be
    \t_{E_n} (\omega) = \min_{A_n\in\aaa (E_n)} \t_{A_n}
    \ee
In words, $E_n$ is a set of patterns, and the
occurrence time of $E_n$ is the volume of
the first cube $B_k$ in which some pattern of $E_n$ can
be found.
A sequence of $\fe_n$-measurable events $E_n$
is called a {\it sequence of rare events} if $\pee (E_n) \to 0$
as $n\to\infty$. For sequences of rare events,
one typically expects so-called {\em exponential laws}, i.e.,
limit theorems of the type
    \be\label{explaw1}
    \pee \left( \t_{E_n} \geq \frac{t}{\pee (E_n)}\right)=
     e^{-\lambda_{E_n} t}+o(1).
    \ee
Equation (\ref{explaw1}) has been proved
for ``high temperature Gibbsian random fields"
and the parameter $\lambda_{E_n}$ is bounded away from zero and infinity.
In the case of patterns, the parameter depends on the self-repetitive
structure of the pattern. For so-called good (meaning that there can be
no fast returns) patterns, we even have that $\lambda_{E_n} =1$. In \cite{acrv}
the exponential law for patterns is
generalized to measurable events $E_n\in\fe_{n}$, provided
a second moment condition is satisfied. This second moment condition ensures that
$1/\pee (E_n)$ is the right time scale for the occurrence
time, i.e., the parameter $\lambda_{E_n}$ is bounded
away from zero and infinity (see Theorem \ref{expthm} below for the precise
formulation).

The relation between maxima and rare events is intuitively obvious:
if a cluster with cardinality bigger than $m$ appears in a cube $B_n$,
then the occurrence time for the appearance of
a cluster bigger than $m$ is not larger than $|B_n|$. More precisely,
define
    \be\label{clustev}
    E_n = \{n\leq |\ce (0)|<\infty\}
    \ee
and define the random variable $\tau_{E_m}$ with values in $\{ (2n+1)^d:n\in\N\}$ by
    \be\label{help}
    \{\tau_{E_m}\leq (2n+1)^d\}= \{ \exists x\in B_n: \theta_x\omega_{B_n}\in E_m \}
    \ee
The random variable $\tau_{E_m}$ is not exactly equal to the occurrence
time $\t_{E_m}$, but we will see that asymptotically $\tau_{E_m}$
and $\t_{E_m}$ have the same distribution (see
Lemma \ref{lem-ttau} below).

The advantage of working with $\tau_{E_m}$ lies in the equality
    \be
    \{\mee_n \geq m\} = \{ \tau_{E_m}\leq (2n+1)^d\}
    \ee
If we assume that the exponential law holds for the occurrence time, then
    \be
    \pee\left(\mee_n \geq u_n(x)\right)\approx 1-
    e^{-\lambda_{E_n} \pee(E_{ u_n (x)}) (2n+1)^d}
    \ee
So if we can choose $u_n (x)$ such that
    \be\label{tune}
    \pee(E_{ u_n (x)})\approx \frac{a_n e^{-x}}{(2n+1)^d}
    \ee
then we obtain (\ref{exval}).
This is the guiding idea of this paper, and the proof
of a result of the type (\ref{exval}) thus relies
on the following three ingredients:

\ben
\item Verification of the validity of the exponential law
for the events $E_n$. For this, we will rely on the techniques developed
in \cite{acrv}, which requires natural mixing conditions
and a second moment estimate, see (cf.\ (\ref{secmcon})).

\item Proof of the existence of the sequence $u_n(x)$ such that
(\ref{tune}) holds.

\item Proof that $\lambda_{E_{u_n(x)}} = 1 + o(1)$.
\een

\section{General results}
\label{sec-genres}
In this section, we introduce the conditions needed and
state the precise form of (\ref{exval}).  We will start
by defining the main conditions in Section 3.1, we will
state the exponential law proved in \cite{acrv} in Section 3.2,
and in Section 3.3, we will state our main results valid under the
formulated conditions.

\subsection{The conditions}
There will be three main
conditions, a non-uniformly exponentially $\phi$-mixing condition,
a finite energy condition, and a condition ensuring that clusters are
subcritical or supercritical.

We first introduce the so-called ``high mixing" condition which is adapted
to the case of Gibbsian random fields.
For $m>0$ define
    \be\label{mixfun}
    \phi(m) = \sup \frac 1{|A_1|}|\ \Pr\left(E_{A_1}|E_{A_2}\right) - \Pr\left(E_{A_1}\right)|,
    \ee
    where the supremum is taken over all finite subsets $A_1, A_2$ of $\Zd$,
with $d(A_1,A_2)\ge m$ and
$E_{A_i}\in\fe_{A_i}$, with $\Pr(E_{A_2})>0$.
Note that this $\phi (m)$ differs from the usual
$\varphi$-mixing function since we divide by the size of the dependence
set of the event $E_{A_1}$. This is natural in
the context of Gibbsian random fields, where the classical $\phi$-mixing mostly
fails (except for the simplest i.i.d.\ case and ad-hoc examples of independent
copies of one-dimensional Gibbs measures).

We are now ready to formulate the non-uniformly exponentially $\phi$-mixing
(NUEM) condition:

\bd[NUEM]
A random field is  {\bf non-uniformly exponentially $\phi$-mixing}
(NUEM)
if there exist  constants $C,c>0$ such that
    \be\label{exponmix}
    \phi(m) \leq  C\exp(-c m)\quad \text{for all}\quad m>0.
    \ee
\ed

Examples of random field satisfying the NUEM condition
are Gibbs measures with exponentially decaying potential
in the Dobrushin uniqueness regime, or local transformations
of such measures. Of course, for site percolation, where we have
independence, we have $\phi =0$.

We next define the finite energy property:

\bd[Finite energy property]
A probability measure $\pee$ has the finite energy property if there exists
$\delta\in (0,1)$ such that
\be\label{fen}
\delta \leq\inf_{\omega\in\Omega} \pee (\omega_x=1|\omega_{\Z^d\setminus \{x \}})
\leq\sup_{\omega\in\Omega} \pee (\omega_x=1|\omega_{\Z^d\setminus \{x \}}) \leq 1-\delta
\ee
\ed

Gibbs measures have the finite energy property
(in particular, it holds of course for independent site percolation, for which
(\ref{fen}) holds with $\delta=1-\delta=p$), but in general
it suffices that there exists a bounded version of
$\log\pee (\si_0 =1|\si_{\{0 \}^c})$.
A direct consequence of
(\ref{fen}) is the existence of $C,C'>0$ such that for any
$\si\in\Omega$, $V\subset\Z^d$,
    \be
    e^{-C|V|}\leq \pee(\omega_V =\si_V ) \leq e^{-C' |V|}.
    \ee
We finally define what it means for a measure to have subcritical
clusters:

\bd[Sub- and supercritical clusters]

\ben
\item[{\rm (i)}]
The probability measure $\pee$ is said to have
subcritical clusters if  $\pee (|\ce (0)| <\infty ) =1$
and if there exists $\zeta, \xi\in (0,\infty)$ such that
    \be\label{subc}
    e^{-\zeta}\leq
    \liminf_{n\to\infty}\frac{\pee(|\ce (0)|\geq n+1)}{\pee(|\ce (0)|\geq n)}
    \leq
    \limsup_{n\to\infty}s\frac{\pee(|\ce (0)|\geq n+1)}{\pee(|\ce (0)|\geq n)}
    \leq e^{-\xi}.
    \ee

\item[{\rm (ii)}]
The probability measure $\pee$ is said to have
supercritical clusters if $\pee (|\ce (0)| =\infty )>0$
and if
    \be\label{supc}
    \lim_{n\rightarrow \infty} \frac{\pee(n+1\leq |\ce (0)|<\infty)}
    {\pee(n\leq |\ce (0)|<\infty)}
    =1.
    \ee
\een
\ed

\subsection{The exponential law}
In order to have the exponential law,
we need that the events $E_n$ are somewhat localized.
More precisely, the non-occurrence of the event
in a big cube can be decomposed as
an intersection of non-occurrence of the event
in a union of small sub-cubes separated by corridors.
Then mixing can be used to factorize the probabilities of on-occurrence
in the sub-cubes, provided the corridors are sufficiently large.
Optimization of this philosophy is the content of the Iteration Lemma
in \cite{acrv}. In our case, the events are not
strictly localized but they can be replaced by local
events, without affecting limit laws. This
is made precise in the following definition:

\bd[Localizability]\label{locl}
\ben
\item[{\rm (i)}] Let $E_n$ be a sequence of events such that
$\pee (E_n)\to 0$. The events are called
local w.r.t. $\pee (E_n)$,
if $E_n\in \fe_{k_n}$ with $k_(2n+1)^d \pee(E_n)^\theta\to 0$
for any $\theta >0$.
\item[{\rm (ii)}] The events $E_n$ of point 1 are called localizable
if there exist events $E'_n$ which are local w.r.t.\
$\pee(E_n)$ such that
$\pee( E_n)-\pee (E'_n)\to 0$ and for any sequence
$u_n\to\infty$
\[
\lim_{n\to\infty}|\pee (\t_{E'_n} \leq u_n)-\pee (\t_{E_n}\leq u_n)|=0
\]
$E'_n$ is then called a local version of $E_n$.
\een
\ed

We will use the following theorem which can
be derived from \cite{acrv}, as we explain below.
\bt[Exponential law]\label{expthm}
Suppose is $\pee$ is finite energy and satisfies the NUEM
condition. Suppose further that $E_n$ are localizable measurable
events such that for some $\delta,\gamma>0$ and all $n\in\N$,
$\pee (E_n) \leq e^{-\gamma n^\delta}$. Assume furthermore that
for any $\alpha>1$
    \be\label{secm}
    \limsup_{n\to\infty}\sum_{0<|x|\leq n^{\alpha}}
    \frac{\pee (E_n\cap \theta_x E_n)}{\pee (E_n)} <\infty
    \ee
then there exists $\La_1,\La_2,c,\rho\in (0,\infty)$, such that
for all $n\in\N$ there exists $\lambda_{E_n}\in [\La_1,\La_2]$ such that
    \be\label{explaw}
    \left|\pee \left(\t_{E_n} > \frac{t}{\lambda_{E_n}\pee (E_n)}\right) - e^{-t}\right|
    \leq \pee (E_n)^\rho e^{-ct}
    \ee
\et

For the ``local version" $E'_n$, the
theorem follows from \cite[Theorem 2.6 and Remark 2.8]{acrv}.
The extension to $E_n$ is  straightforward
from Definition \ref{locl} and is formulated in detail in
\cite[Remark 4.13]{acrv}. Note that there is some notational
difference between the present paper and \cite{acrv}, since
in \cite{acrv}, the occurrence time $\t_{E_n}$ is the {\it width}
of the first cube where $E_n$ occurs, whereas in our setting, it
is the {\it volume}.

Condition (\ref{secm}) is needed to apply Lemma 4.6 in
\cite{acrv}, see also \cite{abadi}. It ensures the existence
of the lower bound $\La_1$ on the parameter $\lambda_{E_n}$
(which is obtained via a second moment estimate for
the number of occurrences). It guarantees further
that the parameter is bounded away from zero which means
that in a cube of volume $\pee (E_n)^{-1}$, the event
$E_n$ happens with a probability bounded away from zero
(uniformly in $n$). This means that $\pee (E_n)^{-1}$ is the
right scale, i.e., a cube with this volume
is such that the event $E_n$ happens with probability
bounded away from zero or one.


The parameter $\lambda_{E_n}$ measures the
``self-repetitive" nature of the event $E_n$, i.e., whether
the event appears typically isolated or in clusters.
See also \cite{abadi} for one-dimensional examples
of $\lambda_{E_n} \not= 1$ and conditions ensuring $\lambda_{E_n} =1$.
For the events $E_n$ of our paper, we will show that $\lambda_{E_n}=1$.


\subsection{Main results}
In our context, Condition (\ref{secm})
is satisfied as soon as for any $\alpha>1$, we have
    \be\label{secmcon}
    \limsup_{n\to\infty}\sum_{0<|x|<n^\alpha}
    \frac{\pee\left(\{n\leq |\ce (x)|<\infty\}\cap \{n\leq |\ce(0)|<\infty\}\right)}
    {\pee(n\leq |\ce(0)|<\infty)} <\infty.
    \ee
The value  of $\alpha$ which we will need later is related to the localization of the
event $|\ce (0)|>n$ to the event $n<|\ce(0)|<n^\alpha$
(see the proof in Section 4 for more details).

Now we can state our main result for the subcritical case:
\bt[Subcritical Gumbel law]\label{mainthm}
Suppose $\pee$ is finite energy, NUEM,
has subcritical clusters and satisfies (\ref{secmcon}).
Then there exists a sequence $u_n\to\infty$, and
a bounded sequence $a_n\in [e^{-\zeta},1]$, such that for $x\geq 0$
    \be
    e^{-a_n e^{-x\zeta}}\leq \pee (\mee_n \leq u_n + x)
    \leq e^{-a_n e^{-x\xi}}.
    \ee
When $x<0$, the upper and lower bound are reversed.

Moreover, if $\xi=\zeta$ , then there exists a constant $\rho>0$ such that
    \be
    |\pee \left(\mee_n \leq u_n + x\right)-e^{-a_n e^{-\zeta x}}|
    \leq \frac{1}{n^\rho}.
    \ee
\et

We now turn to examples where we can apply Theorem \ref{mainthm}.
The following proposition yields a class of non-trivial examples:

\bp\label{exprop}
If $\pee$ is a subcritical Markov measure
satisfying the FKG inequality, then
(\ref{secmcon}) is satisfied.
\ep

This gives the following applications:
\ben
\item Subcritical site percolation $\pee= \pee_p$ where
$\pee_p$ is the Bernoulli measure with $\pee_p (\om_0)=p$
and $p<p_c$.
\item In $d=2$: Ising model at $\beta <\beta_c$. In general
dimension, Ising model at high temperature and/or
high enough magnetic field (see \cite{ghm}).
\een

In very general context we have (\ref{secmcon})
in high enough magnetic field.
The idea is that as soon as for any $V$, and
any $\omega\in\Omega$, the conditional
probabilities $\pee_V (\cdot |\omega_{V^c})$ can
be dominated by a Bernoulli measure
with subcritical clusters, then of course, for any $x\not= 0$,
\eq
\pee \big(|\ce(x)|\geq n\big||\ce (0)|\geq n\big) \leq \pee_p (|\ce (0)| \geq n),
\en
and hence (\ref{secmcon}) is satisfied.

We will now formulate another class of examples.
We say that $\pee$ is dominated by a Bernoulli
measure in the sense of Holley, if for all $\omega\in\Omega$
    \be\label{holdom}
    \pee (\omega_0 = 1|\om_{\{ 0\}^c}) <p.
    \ee
This condition implies that $\pee$ is stochastically dominated by the Bernoulli
measure $\pee_p$.
For measures that are dominated by a subcritical Bernoulli
measure, our results also apply:

\bp
Let $p_c$ denote the critical value for Bernoulli
site percolation.
If (\ref{holdom}) is satisfied for
some $p<p_c$, then (\ref{secmcon}) holds
true.
\ep

This proposition can be applied to Gibbs measures
such that the potential has a Dobrushin norm
which is small enough (to guarantee mixing condition),
with magnetic field high enough such that (\ref{holdom})
holds, see \cite{ghm} for  more details.

Our last theorem applies for independent supercritical
site percolation. Recall (\ref{meenzbsup}).
Then we have the following result for supercritical site percolation:

\bt[Supercritical Gumbel law]\label{super}
Let $p>p_c$. There exists a constant $a>0$, a sequence
$a_n\in (a,1]$ and a sequence $u_n(x)$ such that
$u_n(x)\uparrow\infty$ for all $x$ as $n\uparrow\infty$,
such that for all $x$
    \be
    \pee_p (\mee_n\leq u_n(x)) = e^{-a_n e^{-x}}+o(1).
    \ee
If $\pee$ has finite supercritical clusters, then $a_n=1+o(1)$.
\et

\section{Proofs}
\label{sec-proofs}
In this section, we prove the main results stated in Sections
\ref{sec-intro} and \ref{sec-genres}.

\subsection{Preparations}
In this section, we state some general results for non-critical
clusters. In Lemma \ref{clust}, we first identify the law of
$|\ce (0)|$ for non-critical clusters in terms of the law
of $|\cee (0)|$. In Proposition
\ref{prop-ratiops} and Lemma \ref{seqlem}, we investigate the
cluster size distribution in more detail.

\bl[The law of $|\ce(0)|$]\label{clust} Suppose that
$\mathbb{E}_p\left(|\cee(0)| I[|\cee (0)| <\infty]\right)<\infty$. Then, for $n\geq 1$,
    \be
    \pee (|\cee (0)| =n)= n \pee (|\ce (0)| =n).
    \ee
\el

\bpr
We start with the subcritical case.
Let $V_k$ be a sequence of volumes such
that $V_k\to\Z^d$ and $|\partial V_k|/|V_k| \to 0$
as $k\to\infty$. For a cluster $\ce$, we
denote by $le(\ce)$  the left endpoint of $\ce$,
in particular $le(\ce (x))=x$ by definition.
Then, for $n\geq 1$,
    \beq
    |V_k|n\pee (|\ce(0)| =n)
    &=& \sum_{x\in V_k} n \pee (|\ce (x)| =n)
    \nonumber\\
    &=&
    \sum_{x\in V_k} \sum_{y} \pee (|\ce (x)| =n, y\in \ce(x))
    \nonumber\\
    &=&
    \sum_{x\in V_k}\sum_{y} \pee \left(|\cee (y)| =n
    , x= le(\cee(y))\right)
    \nonumber\\
    &=&
    \sum_{y}\pee\left( |\cee(y)| =n, le (\cee (y))\in V_k\right)
    \nonumber\\
    &=& |V_k| \pee ( |\cee (0)|=n) + \sum_{y\not \in V_k}
    \pee\left( |\cee(y)| =n, le (\cee (y))\in V_k\right)\nonumber\\
    &&\quad -
    \sum_{y\in V_k }\pee\left( |\cee(y)| =n, le (\cee (y))\not \in V_k\right).
    \eeq
We claim that the last two terms are $O (|\partial V_k|)$. Indeed, the
first sum is equal to
    \be
    \sum_{y\not\in V_k}\sum_{x\in V_k}  \pee (|\ce(x)| =n, y\in\ce (x))
    \leq \E \Big|\bigcup_{x\in V_k: \ce(x)\cap V_k^c\neq \varnothing} \ce(x)\Big|
    \leq \sum_{z\in \partial V_k} \E|\cee(z)|=O (|\partial V_k|),
    \ee
where the last step follows because $\pee$ has subcritical clusters,
and hence $\pee ( |\cee (x) | \geq n) \leq e^{-\zeta n}$ by
\cite[Theorem (6.78)]{Grim99}, so that $\E|\cee(z)|<\infty$
for all $z$. The second sum
is bounded similarly as
    \be
    \sum_{y\in V_k}\sum_{x\not\in V_k}  \pee (|\ce(x)| =n, y\in \ce (x))
    \leq \E \Big|\bigcup_{x\not\in V_k: \ce(x)\cap V_k\neq \varnothing} \ce(x)\Big|
    \leq \sum_{z\in \partial V^c_n} \E|\cee(z)|=O (|\partial V_k|).
    \ee
Therefore, we obtain that
    \eq
    |V_k|n\pee (|\ce(0)| =n) = |V_k| \pee ( |\cee (0)|=n) + O (|\partial V_k|).
    \en
Divide by $|\partial V_k|$ and let $k\rightarrow \infty$ to arrive at the claim.

The proof for the supercritical case is similar and based on the estimate
$\pee (n\leq |\cee (x)|<\infty)\leq e^{-\gamma n^\delta}$ from which we
conclude $\E (|\cee (z)| I [|\cee (z)| <\infty]) <\infty$.
\epr

Before we formulate our next proposition, we remark that the cluster
$\ce(0)$ is finite with probability one, since
it has $0$ as its left endpoint. Therefore, we have
$\pee (n\leq |\ce (0)|<\infty)=\pee (|\ce (0)| \geq n) $
in the supercritical case, and we can drop the restriction
that the cluster is finite in the notation. Naturally, we also drop
this restriction in the subcritical case.

\bp[Lower bound on the cluster tail] \label{prop-ratiops} If
$\pee$ is finite energy, then there exists a $\zeta >0$ such that
\be\label{511} \liminf_{n\to\infty} \frac{\pee (|\ce(0)|\geq
n+1)}{\pee (|\ce(0)|\geq n)} \geq e^{-\zeta}. \ee \ep

\bpr
We start with the subcritical case.
We abbreviate $X= |\ce (0)|$. Define $A_n = \{X\geq n\}$.
The ratio we are interested in can be written as
    \be
    \frac{\pee (X\geq n+1)}{\pee (X\geq n)}
    =
    \frac{\int_{A_{n+1}} d\pee}{\int_{A_{n+1}} d\pee +
    \int_{A_n\setminus A_{n+1}} d\pee}
    \ee
For $\omega\in\Omega$, we denote by $\omega^x$ the configuration
obtained by flipping at $x$, i.e.,
$\omega^x (y) = \omega (y)$ for $y\not= x$ and
$\omega^x (x) = 1-\omega (x)$.
Then, for any $\om\in A_n\setminus A_{n+1}$ there exists $x= x_\om\in\Z^d$
such that $\om^x\in A_{n+1}$, which gives
    \be\label{krak}
    \frac{\int_{A_{n+1}} d\pee}{\int_{A_{n+1}} d\pee +
    \int_{A_n\setminus A_{n+1}} d\pee}
    \geq
    \frac{\int_{A_{n+1}} d\pee}{\int_{A_{n+1}} d\pee +
    \int_{A_{n+1}} d\pee^{x_\om}}
    \ee
where for $x\in\Z^d$, $\pee^{x}$ denotes the image measure
under the transformation $\Psi_x: \om\mapsto\om^{x}$.
The finite energy property implies
that
    \be\label{lap}
    C= \sup_x\|\frac{d\pee^x}{d\pee} \|_{L^\infty (\pee) }
    <\infty.
    \ee
Therefore, from (\ref{krak}), we obtain the lower bound
    \[
    \frac{\pee (X\geq n+1)}{\pee (X\geq n)}
    \geq
    \frac{1}{1+C},
    \]
which is (\ref{511}).

To deal with the supercritical case, let $\omega\in A_n\setminus A_{n+1}$.
Flipping one occupation variable at the exterior boundary of $\cee (0)$
can lead to an infinite cluster, so that $\omega^x\not\in A_{n+1}$. However,
we can make one site occupied and all its neighbors which do not belong
to $\ce (0)$ vacant. This leads to a configuration in $A_{n+1} =
\{n+1\leq |\ce (0)|<\infty\}$. Since this transformation
$T_x$ of $\omega$ is still local (it affects only one site
$x=x_\omega$ and possibly some neighbors), the same
argument applies where now, using the finite energy
property we replace $C$ of (\ref{lap}) by
    \be
    C' = \sup_x\|\frac{d\pee\circ{T_x}}{d\pee} \|_{L^\infty (\pee) }
        <\infty.
    \ee

\epr

\bl[Existence of $v_n$]\label{seqlem} There exists
$v_n\uparrow\infty$ and a sequence $b_n$ satisfying $0<b_n\leq 1$
such that
    \be
    \pee ( |\ce (0)| > v_n ) = \frac{b_n}{n}.
    \ee
It $\pee$ has subcritical clusters, then
$\liminf_{n\rightarrow \infty} b_n\geq e^{-\zeta}$,
while if $\pee$ has finite supercritical clusters,
then $b_n=1+o(1)$.
\el

\bpr
We abbreviate $X= |\ce (0)|$, and recall that $X<\infty$ both
above and below criticality.  We define
    \beq
    v_n^+ &=& \inf \{ x\in \N: \pee (X \geq x) \leq \frac{1}{n}\},
    \nonumber\\
v_n^- &= & \sup \{ x\in \N: \pee (X \geq x) \geq \frac{1}{n}\}.
    \eeq
Then $v_n^+ = v_n^- +1$ or $v_n^+ = v_n^-$. Put $v_n = v_n^+$.
By definition
    \[
    \pee (X \geq v_n) \leq \frac1n,
    \]
so that $b_n\leq 1$.

Moreover,
    \beq
    n\pee (X\geq v_n) &=& n \pee (X\geq v_n^-)
    \frac{ \pee(X\geq v_n^+)}{\pee(X\geq v_n^-)}
    \nonumber\\
    &\geq &
    \frac{ \pee(X\geq v_n^- +1)}{\pee(X\geq v_n^-)}.
    \eeq
We note that $v_n\rightarrow \infty$ when $n\rightarrow \infty$.
Therefore, if $\pee$ has subcritical or supercritical clusters,
    \be
    \liminf_{n\to\infty} n \pee (X\geq v_n) \geq
    \liminf_{n\to\infty}\frac{ \pee (X\geq n+1)}{\pee (X\geq n)} = e^{-\zeta}.
    \ee
Thus, we obtain $\liminf_{n\to\infty}{b_n} \geq e^{-\zeta}$.
On the other hand, when $\pee$ has finite supercritical clusters,
    \be
    1\geq \lim_{n\to\infty} n \pee (X\geq v_n) \geq
    \lim_{n\to\infty}\frac{ \pee (X\geq n+1)}{\pee (X\geq n)} = 1.
    \ee
Thus, we obtain $\lim_{n\to\infty}{b_n} =1$.
\epr

We next verify that the events
    \eq
    E_n  = \{|\ce (0)| \geq n\}
    \en
are localizable. This is the content of the next lemma.

\bl[Localizability of $E_n$]
\label{lem-locsub}
The events $E_n  =\{|\ce (0)| \geq n\}$ are
localizable, and their local versions
can be chosen as
    \eq
    E'_n = \{n\leq |\ce (0)|< n^\theta\}
    \en
for some $\theta\in (1,\infty)$ with $k_n=n^{\theta}$ in
Definition \ref{locl}.
\el

\bpr This is an obvious consequence  of the estimates that there exist positive
$\zeta=\zeta(p,d)$ and $\eta=\eta(p,d)$ such that $\pee (|\cee
(0)|>n ) \leq e^{-\zeta n}$ in the subcritical case and
$\pee(n\leq |\cee(0)|<\infty ) \leq e^{-\eta n^{\frac{d-1}{d}}}$
in the supercritical case, for some $\zeta,\eta>0$. See
\cite[Theorems (6.78) and (8.61)]{Grim99} for these
estimates in the context of percolation.
\epr

We finish this section with a lemma showing the asymptotic equivalence
of $\tau_{E_n}$ introduced in (\ref{help}) and the occurrence time
$\t_{E_n}$.


More precisely, we have the following lemma:

    \bl[Occurence times]
    \label{lem-ttau}
    Let $m=m_n\uparrow\infty$ be such that $m_n n^{\epsi-1}$ converges to zero as $n\to\infty$ for
    some $\epsi\in (0,1)$, and
    such that $\pee (E_{m_n})\leq  n^{-d +\epsi}$. Then
    \be
    \pee (\tau_{E_{m_n}}\leq (2n+1)^d )=\pee (\t_{E_{m_n}}\leq (2n+1)^d )+o(1).
    \ee
    \el

\bpr
First we remark that
    \be
    \{\t_{E_m}\leq (2n+1)^d\} \subset \{ \tau_{E_m}\leq (2n+1)^d\},
    \ee
and
    \be
    \{\tau_{E_m}\leq (2n+1)^d\} \setminus \{ \t_{E_m}\leq (2n+1)^d\}\subset
    \{ \exists x\in B_n: x+B_m\not\subset B_n: |\ce(x)|>m\}.
    \ee
We estimate
    \beq
    \pee\left(\{ \exists x\in B_n: x+B_{m_n}\not\subset B_n: |\ce(x)|>m_n\}\right)
    & \leq &
    \pee (|\ce(0)|\geq m_n) |\{x\in B_n:x+B_{m_n}\not\subset B_n\}|\nonumber\\
    & \leq & \frac{m_n n^{d-1}}{n^{d-\epsi}}=m_n n^{\epsi-1}.
    \eeq
This converges to zero as $n\to\infty$ by the assumption on $m_n$.
\epr

\subsection{Maximal subcritical clusters}

In this section, we prove Theorems \ref{thm-Gumbel} and \ref{mainthm}.
We  study the tails of the cluster size distribution, subject to
(\ref{subc}). The main result is the following lemma:

\bl[Identification $u_n(x)$] \label{lem-taildistrsub} Suppose
$\pee$ has finite energy and has subcritical clusters, then there
exists a sequence $a_n$ with $0<a_n\leq 1$ such that
$\liminf_{n\rightarrow\infty} a_n \geq e^{-\zeta}$, such that for
all $x >0$ and for all $n\geq 1,$
    \be
    \frac{a_n}{(2n+1)^d}e^{-\zeta x}(1+o(1))\leq \pee (|\ce (0)|\geq
    u_n +x)\leq  \frac{a_n}{(2n+1)^d}e^{-\xi x}(1+ o(1)).
    \ee
For $x<0$ the same inequality holds with $\zeta$ and $\xi$ interchanged.
\el

\bpr
Let $x>0$. We again abbreviate $X = |\ce (0)|$.
    \be
    \pee (X \geq u_n +x)= \prod_{i=1}^{x}
    \frac{\pee (X \geq u_n +i)}{\pee (X \geq u_n +i-1)}.
    \ee
Hence, with the choice of $u_n= v_{(2n+1)^d}$ where $v_n$ is
as in Lemma \ref{seqlem},
    \be
    \liminf_{n \to\infty} n \pee (X \geq u_n +x)
    \geq \liminf_{n\to\infty} n \pee (X\geq u_n) e^{-\zeta x},
    \ee
and, using Lemma \ref{seqlem} again, for any $x\in \N$ fixed,
    \be
    \limsup_{n \to\infty} n \pee (X \geq u_n +x)
    \leq \left(\limsup_{n\to\infty}
    \frac{\pee (X\geq n+1)}{\pee (X\geq n)}\right)^x
    = e^{-\xi x}.
    \ee
This proves the claim for $x>0$. The proof for $x<0$ is similar.
\epr


We now verify Condition (\ref{secm}) for FKG measures.\\
\\
{\bf Proof of Proposition \ref{exprop}.} We have to prove that
for any $\alpha>0$,
    \be
    \label{goal}
    \limsup_{n\to\infty}\sum_{0<|x|\leq n^{\alpha}}
    \frac{\pee (|\ce(0)|\geq n, |\ce (x)|\geq n)}
    {\pee (|\ce(0)|\geq n)} <\infty.
    \ee
In fact, we will show that the right-hand side of (\ref{goal}) converges to 0 when
$n\rightarrow \infty$.

We denote by $\pee_\la^\eta(\omega_\la)$ the conditional probability
to find $\omega$ inside $\la$, given $\eta$ outside $\la$. For a Markov random
field, the dependence on $\eta$ is only through the boundary of $\la$, i.e.,
    \eq
    \pee_\la^\eta(\omega_\la)=\pee_{\la}^{\eta_{\partial\la}}(\omega_\la),
    \en
where $\partial\la$ denotes the exterior boundary of $\la$, i.e., the set of
those sites not belonging to $\la$ which have at least one neighbor inside $\la$.
Thus, we can think of $\eta$ as describing the boundary condition.
By the FKG-property, we have that if $\eta\leq \zeta$ for $\eta,
\zeta \in \{0,1\}^{\Z^d}$, then
    \be\label{fkg}
    \pee_\la^\eta\leq \pee_\la^\zeta
    \ee
Moreover, by definition of the clusters $\ce (x)$, $\ce(0)\cap \ce(x)= \emptyset$
for $x\not=0$. Therefore, we can write, for $x\neq 0$,
    \beq\label{kwats}
    &&\pee (|\ce(0)|\geq n, |\ce (x)|\geq n)\nonumber\\
    &&\qquad =
    \sum_{A: |A|\geq n, le(A)=0} \pee (\ce(0) =A, |\ce (x)|
    \geq n, \ce (x)\cap A=\emptyset)
    \nonumber\\
    &&\qquad =
    \sum_{A: |A|\geq n, le(A)=0}
    \pee (|\ce (x)| \geq n, \ce (x)\cap A=\emptyset|\om_A=1,
    \om_{\partial A} =0)\pee (\ce(0)=A)
    \nonumber\\
    &&\qquad \leq \sum_{A: |A|\geq n, le(A)=0}
    \pee (\ce(0)=A)\pee^{0_{\partial A}}_{\Z^d\setminus \bar{A}} (|\ce (x)|\geq n),
    \eeq
where in the last step we have used the Markov property,
with the notation $\bar{A}= A\cup\partial A$.
Using (\ref{fkg}), we thus arrive at
    \be
    \label{FKGapplic}
    \pee^{0_{\partial A}}_{\Z^d\setminus \bar{A}} (n\leq |\ce (x)|\leq n^{\theta})
    \leq \pee_{\Z^d\setminus \bar{A}} (|\cee (x)|\geq n)
    \leq \pee(|\cee (x)|\geq n).
    \ee
Equation (\ref{FKGapplic}) combined with (\ref{kwats}) leads to the correlation inequality
    \be
    \pee(|\ce (0)|\geq n, |\ce (x)|\geq n)
    \leq \pee (|\ce (0)|\geq n)\pee (|\cee (0)|\geq n).
    \label{FKG-applic}
    \ee
Therefore,
    \be
    \sum_{0<|x|\leq n^{\alpha}}
    \frac{\pee (|\ce(0)|\geq n, |\ce (x)|\geq n)}
    {\pee (|\ce(0)|\geq n)} \leq
    (2n^{\alpha}+1)^d \pee (|\cee (0)|\geq n) \rightarrow 0,
    \ee
because the decay of the probability $\pee (|\cee (0)|\geq n )$
is faster than $\frac{1}{n^\beta}$ for any $\beta >0$.
\QED

\bp[The subcritical intensity is one]\label{prop-lambsub} For
$u_n$ as in Lemma \ref{seqlem} and for every $x$ bounded, there
exists a $\beta>0$ such that
    \eq
    1-\pee(E_{u_n+x})^{\beta}\leq \lambda_{E_{u_n+x}}\leq 1.
    \en
\ep

\bpr We will first identify $\lambda_{E_{u_n+x}}$. We use \cite[(2.6)]{acrv}, which states
that
    \eq
    \label{lambdadef}
    \lambda_{E} = -\frac{\log \pee\big(\t_E> f_E\big)}{f_{E}\pee(E)},
    \en
where, for some $\gamma\in (0,1)$
    \eq
    f_E = \lfloor \pee(E)^{-\gamma}\rfloor.
    \en
We will show that $\pee\big(\t_E\leq f_E\big)$ is quite small
(as proved in the sequel), so that we can approximate
    \eq
    -\log\pee\big(\t_E> f_E\big)=\pee\big(\t_E\leq f_E\big)
    + O\Big(\pee\big(\t_E\leq f_E\big)^2\Big).
    \en
Therefore,
    \eq
    \lambda_{E} = \frac{\pee\big(\t_E\leq f_E\big)}{f_{E}\pee(E)} +o(1).
    \label{lambdaasy}
    \en

We will proceed by computing $\pee\big(\t_E\leq f_E\big)$. To do so, we write
    \eq
    \pee\big(\t_{E_{u_n+x}}\leq f_{E_{u_n+x}}\big) = \pee\big(\bigcup_{y\in B_{m_{n,x}}}
    \{|\ce(y)|\geq u_n+x\}),
    \en\
where we abbreviate $m_{n,x}=f_{E_{u_n+x}}^{1/d}$.
By Boole's inequality,
    \eq
    \label{boole}
    \pee\big(\t_{E_{u_n+x}}\leq f_{E_{u_n+x}}\big) \leq \sum_{y\in B_{m_{n,x}}} \pee(|\ce(y)|\geq u_n+x)
    =f_{E_{u_n+x}}\pee(E_{u_n+x}).
    \en
Thus,
    \eq
    \lambda_{E_{u_n+x}}\leq 1.
    \en
For the lower bound, use
    \beq
    \label{incl-excl}
    \pee\big(\t_{E_{u_n+x}}\leq f_{E_{u_n+x}}\big) &\geq& \sum_{y\in B_{m_{n,x}}}\pee(|\ce(y)|\geq u_n+x)\\
    &&\qquad -\shift\sum_{y, z\in B_{m_{n,x}}: y\neq z}\shift
    \pee(|\ce(y)|\geq u_n+x, |\ce(z)|\geq u_n+x).\nonumber
    \eeq
The first term is identical to the first term in the upper bound, and we need
to bound the second term only. For this, we use (\ref{FKG-applic}), and thus
obtain
    \be
    \pee\big(\t_{E_{u_n+x}}\leq f_{E_{u_n+x}}\big)
   \geq f_{E_{u_n+x}}\pee(E_{u_n+x})
    -f_{E_{u_n+x}}^2\pee(E_{u_n+x}) \pee(|\cee(0)|\geq u_n+x).
    \ee
Thus,
    \eq
    \lambda_{E_{u_n+x}}\geq 1-f_{E_{u_n+x}}\pee(|\cee(0)|\geq u_n+x) \geq 1-\pee(E_{u_n+x})^{\beta}
    \en
for some $\beta>0$.
\epr

We finally identify the sequence
$u_n$ under the hypothesis of a ``classical" subcritical
cluster tail behavior in Proposition \ref{sublem}, and
under the hypothesis of a ``classical" supercritical
cluster tail behavior in Proposition \ref{suplem}.

\bp[Identification $u_n(x)$ for classical subcritical
tails]\label{sublem}~\\
Suppose that there exists $\alpha\in \R,
\zeta>0$ and $0<C<\infty$, such that
    \be
    \pee (|\ce (0)|\geq n )= Cn^\alpha e^{-\zeta n}[1+o(1)].
    \ee
Then
    \be
    u_n =  \Big\lfloor \frac{\log n}{\zeta} + \frac{\alpha \log\log n}{\zeta}\Big\rfloor.
    \ee
\ep

\bpr
This is a simple computation, using the definition
of $u_n$ introduced in the proof of Lemma \ref{seqlem}.
\epr

\vskip0.5cm

\noindent
{\bf Proof of Theorem \ref{mainthm} and Theorem \ref{thm-Gumbel}.}
We first finish the proof of Theorem \ref{mainthm}.
We first use the equality
    \eq
    \label{Mntaueq}
    \{\mee_n \geq m\} = \{ \tau_{E_m}\leq (2n+1)^d\}.
    \en
Then we use Lemma \ref{lem-ttau} to obtain that as long as
$\pee (E_{u_n(x)})\leq n^{-d+\epsi}$, we have
    \eq
    \label{tothet}
    \pee (\mee_n\geq u_n(x))=\pee(\t_{E_{u_n(x)}}\leq (2n+1)^d)+o(1).
    \en
We wish to apply Theorem \ref{explaw}, and will first check that the
conditions are fulfilled. We note from Lemma \ref{lem-locsub} that the
events $E_{u_n(x)}$ are localizable with local versions
$E_{u_n(x)}'$. Furthermore, from Proposition \ref{exprop},
it follows that Condition (\ref{secm}) is fulfilled for $E_{u_n(x)}$.
Therefore, we may apply Theorem \ref{explaw}.

We choose $u_n(x)=u_n+x$ as in Lemma \ref{unx}, and the event
$E_{u_n+x}$ as before. Note that for this $u_n(x)$, we indeed
have that for every $x$ fixed,
    \eq
    \pee (E_{u_n+x})=\frac{e^{-x}}{(2n+1)^d} a_n \leq n^{-d+\epsi},
    \en
so that we can use (\ref{tothet}).

Assume that $x\geq 0$. For $x<0$ some inequalities
reverse sign. Then we apply Theorem \ref{explaw} to obtain:
    \be
    \pee (\mee_n\geq u_n+x) = \pee(\t_{E_{u_n+x}}\leq (2n+1)^d)+o(1)=
    1-\exp\left(-\lambda_{E_{u_n+x}}(2n+1)^d\pee (E_{u_n+x})\right)+o(1).
    \ee
We need to investigate the exponent. By Lemma \ref{lem-taildistrsub}, we have that
    \eq
    \frac{a_n}{(2n+1)^d}e^{-\zeta x}\leq \pee (E_{u_n+x})
    \leq \frac{a_n}{(2n+1)^d}e^{-\xi x},
    \en
and this inequality is reversed for $x<0$. By Proposition \ref{prop-lambsub},
we have that
    \eq
    \lambda_{E_{u_n+x}}=1+o(1).
    \en
Therefore, for any $x\in\N$,
    \be
    1- \exp(-a_n e^{-\xi x}) + o(1)\leq
    \pee (\mee_n\geq u_n+x)\leq 1- \exp(-a_n e^{-\xi x}) + o(1).
    \ee
This completes the proof of Theorem \ref{mainthm}. When $\zeta=\xi$,
the statement in Theorem \ref{thm-Gumbel} is a direct consequence of
Theorem \ref{mainthm}, combined with Lemma \ref{lem-ttau}.

\QED

\noindent {\bf Remark.}
The examples mentioned in Section 1.3.2 fit into the context of
Theorem \ref{mainthm}. Indeed, for the Ising model, the inequality
(\ref{subc}) is verified above the critical temperature
in $d=2$ and at high enough temperature in any dimension. The mixing
condition (\ref{exponmix}) is verified at high temperature
in the Dobrushin uniqueness regime, and in $d=2$ above the critical
temperature, by complete analyticity.
For general Gibbs measures with a potential with a finite Dobrushin norm,
one can choose the magnetic field high enough such that the Dobrushin
uniqueness condition and hence condition
(\ref{exponmix}) is satisfied (see e.g. \cite{Geo}), and such that
(\ref{subc}) follows from a domination with Bernoulli measures
(see \cite{ghm}).



\subsection{Maximal supercritical clusters}
In this subsection we prove Theorems \ref{super} and \ref{thm-Gumbelsup}.

In the following proposition we show that we can still find a
sequence $u_n(x),$ but not necessarily of the form
$u_n +x$, if we omit the subcriticality condition.
This will be useful when we study the supercritical
percolation clusters.

\bl[Existence of $u_n(x)$]\label{unx} Suppose $\pee$ is finite
energy, NUEM and $\mathbb{E}(|\cee(0)| I [|\cee(0)| <\infty]) <\infty$. Then there exists a
function $u_n(x)$ such that
    \be
    \pee (|\ce(0)| \geq u_n(x)) = \frac{e^{-x}}{(2n+1)^d} a_n,
    \ee
where $a_n$ is a bounded sequence (not depending on $x$). Furthermore, if
$\pee$ has finite supercritical clusters, then $a_n=1+o(1)$.
\el
\bpr
Since $\mathbb{E}(|\cee(0)| I [|\cee(0)| <\infty]) <\infty$, we can use Lemma
\ref{clust}. As in the proof of Lemma \ref{seqlem}, we define
    \beq
    v^+_n(x) &=& \inf \{ k: \pee (X\ge k) \geq\frac{e^{-x}}{n}\},
    \nonumber\\
    v^-_n(x) &=& \sup \{ k: \pee (X\ge k) \leq
    \frac{e^{-x}}{n} \}.
    \eeq
Then we can repeat the proof of Lemma
\ref{seqlem}, and use (\ref{supc}) to conclude that
$a_n = 1 + o(1)$. We can then choose $u_n = v_{(2n+1)^d}$.
\epr

We continue with the following proposition which
will guarantee Condition (\ref{secm}) for
finite supercritical clusters.

%

%
%
%
%
%
%

    \bp[Supercritical second moment condition] For every $\alpha>1$
    \label{prop-secmomsup}
        \be\label{supersec}
        \limsup_{n}\sum_{0<|x|<n^{\alpha}} \frac{\pee \left(E_n\cap \theta_x E_n
        \right)}{\pee\left(E_n\right)}=0.
        \ee
    \ep

\bpr We rewrite
    \eq
    \pee \left(E_n\cap \theta_x E_n
        \right)
        =
        \pee \left(E_n\cap \theta_x E_n\cap \{x\nc \partial
        B_{2n^{\alpha}}\}\right)+\pee \left(E_n\cap \theta_x E_n\cap \{x\conn
        \partial B_{2n^{\alpha}}\}\right).
        \en
The second term is simple, since the event is contained in the
probability that $\{|\cee(x)|>n^\alpha\}$. We bound its contribution
to the left-hand side of (\ref{supersec}) by
    \eq
    (n^{\alpha}+1)^d \frac{\pee \left(E_{n^{\alpha}}\right)}
    {\pee\left(E_n\right)},
    \en
which is an error for any $\alpha>1$.

For $\Gamma\subset\Z^d$, we denote
$\overline{\Gamma}=\Gamma\cup\partial_e\Gamma$. Then, we compute
    \beq
    &&\pee \left(E_n\cap \theta_x E_n\cap \{x\nc \partial
    B_{2n^{\alpha}}\}\right)\\
    &&\qquad =
    \pee \left(E_n\cap \{|\ce(x)|\geq n\}\cap \{x\nc \partial
    B_{2n^{\alpha}}\}\right)\nonumber\\
    &&\qquad =
    \sum_{\Gamma\in {\cal G}_n(0)}
    \pee (\ce (0)=\Gamma)
    \pee\left(|\ce(x)|\geq n,x\nc \partial
    B_{2n^{\alpha}}|\ce (0) = \Gamma\right)
    \nonumber\\
    &&\qquad =
    \sum_{\Gamma\in {\cal G}_n(0)}
    \pee (\ce (0)=\Gamma)
    \pee_{\Z^d\setminus\overline{\Gamma}}\left(|\ce(x)|\geq n,x\nc \partial
    B_{2n^{\alpha}}|\ce (0) = \Gamma\right),
    \nonumber
    \eeq
where
    \eq
    {\cal G}_n(0)=\{\Gamma: 0=le(\Gamma), n\leq |\Gamma|<\infty\},
    \en
and where $\pee_{\Z^d\setminus\overline{\Gamma}}$ is the
conditional measure given that all sites in $\overline{\Gamma}$ are
vacant. We further define
    \be
    \cere_x =
    \cup_{y\in\partial (x+ B_{2n^{\alpha}})} \left(\cee (y) \cap (x+B_{2n^{\alpha}})\right).
    \ee
Then we can further condition on $\cere_x$:
    \beq
    &&\sum_{\Gamma\in {\cal G}_n(0)}
    \pee (\ce (0)=\Gamma)
    \pee_{\Z^d\setminus\overline{\Gamma}}\left(|\ce(x)|\geq n,x\nc \partial
    B_{2n^{\alpha}}|\ce (0)=\Gamma\right)
    \nonumber\\
    &= &
    \sum_{\Gamma\in {\cal G}_n(0)}\sum_{CR\ldots}
    \pee(\ce (0)=\Gamma)\pee (\cere_x =CR)
    \pee_{\Z^d\setminus (\overline{\Gamma}\cup\overline{CR})} (|\ce (x)|\geq n)\nonumber
    \eeq
where we abbreviated the conditions the set $CR$ has to satisfy by
$\ldots$. We can then proceed, using the FKG inequality
    \beq
    && \sum_{\Gamma\in {\cal G}_n(0)}\sum_{CR\ldots}
    \pee(\ce (0)=\Gamma)\pee (\cere_x =CR)
    \pee_{\Z^d\setminus (\overline{\Gamma}\cup\overline{CR})} (|\ce (x)|\geq n)
    \nonumber\\
    &&\qquad\leq
    \sum_{\Gamma\in {\cal G}_n(0)}\sum_{CR\ldots}
    \pee(\ce (0)=\Gamma)\pee (\cere_x =CR)
    \pee_{\Z^d\setminus (\overline{CR})} (|\cee (x)|\geq n)
    \nonumber\\
    &&\qquad \leq
    \sum_{\Gamma\in {\cal G}_n(0)} \pee (\ce (0)= \Gamma)
    \pee(n\leq |\cee (x)|<\infty)\nonumber\\
    &&\qquad =\pee \left(E_n\right)\pee(n\leq |\cee (0)|<\infty).
    \label{FKGbd}
    \eeq
By \cite[Theorem (8.65)]{Grim99}, there exists $\eta=\eta(p,d)>0$ such that
    \be
    \label{clusbd}
    e^{-\gamma n^{\frac{d-1}{d}}} \leq \pee \left(E_n\right)\leq
    \pee \left(n\leq |\cee (0)| <\infty\right)\leq e^{-\eta n^{\frac{d-1}{d}}}.
    \ee
From (\ref{FKGbd}) and (\ref{clusbd}), we conclude that there
exists $\delta>0$ such that
    \beq
    \sum_{0<|x|<n^{\alpha}} \frac{\pee \left(E_n\cap \theta_x E_n
    \right)}{\pee\left(E_n\right)}
    &\leq& (2n^{\alpha}+1)^d \pee(n\leq |\cee (0)|<\infty)
    +(2n^{\alpha}+1)^d\frac{\pee(E_{n^\alpha})}{\pee(E_n)}
    \nonumber\\
    &\leq& (2n^{\alpha}+1)^d e^{-\delta
    n^{\frac{d-1}{d}}},
    \eeq
and thus (\ref{supersec}) follows.\\

\epr

\bp[The supercritical intensity is one]
\label{prop-lambsup} For
$u_n(x)$ as in Lemma \ref{seqlem} and for every $x$ bounded, there
exists a $\beta>0$ such that
    \eq
    1-\pee(E_{u_n(x)})^{\beta}\leq \lambda_{E_{u_n(x)}}\leq 1.
    \en
\ep

\bpr We follow the proof of Proposition \ref{prop-lambsub}.
We will first identify $\lambda_{E_{u_n(x)}}$. Recall (\ref{lambdadef})
and (\ref{lambdaasy}). The upper bound in (\ref{boole}) applies verbatim.

For the lower bound,  use
    \beq
    \pee\big(\t_{E_{u_n(x)}}\leq f_{E_{u_n(x)}}\big) &\geq&
    \sum_{y\in B_{m_{n,x}}}\pee(u_n(x)\leq |\ce(y)|<\infty)\\
    &&\qquad -\shift\sum_{y, z\in B_{m_{n,x}}: y\neq z}\shift
    \pee(u_n(x)\leq |\ce(y)|<\infty, u_n(x)\leq |\ce(z)|<\infty),\nonumber
    \eeq
where now $m_{n,x}=f_{E_{u_n(x)}}^{1/d}$.
The first term is identical to the first term in the upper bound, and we need
to bound the second term only. In order to do so, we derive a similar
bound as in (\ref{FKG-applic}), which was instrumental in the proof of
Proposition \ref{prop-lambsub}.


Write
    \beq
    \shift&&\pee(u_n(x)\leq |\ce(y)|<\infty, u_n(x)\leq |\ce(z)|<\infty)\\
    &&\qquad =\pee(|\ce(y)|\leq u_n(x), |\ce(z)|\leq u_n(x), y,z \nc \partial B_{n})
    \nonumber\\
    &&\qquad\qquad+\pee\big(\{u_n(x)\leq |\ce(y)|<\infty\}\cap \{u_n(x)\leq |\ce(z)|<\infty\}\cap
    \big(\{y \conn \partial B_{n}\}\cup \{z\conn \partial B_{n}\}\big)\big).
    \nonumber
    \eeq
The first is bounded by a similar argument as in (\ref{FKGbd}) by
    \beq
    \shift&&\pee(|\ce(y)|\geq u_n(x), |\ce(z)|\geq u_n(x), y,z \nc \partial B_{n})
    \label{yzbd}\\
    &&\qquad \leq \pee(|\ce(y)|\geq u_n(x), y\nc \partial B_{n})
    \pee(|\cee(z)|\geq u_n(x), z\nc \partial B_{n})\nonumber\\
    &&\qquad \leq
    \pee(E_{u_n(x)})\pee(u_n(x)\leq |\cee (0)|<\infty).\nonumber
    \eeq
Using that $u_n(x)\leq n/2$, the second event in (\ref{yzbd}) is bounded by
    \[\pee(E_{n-u_n(x)})\leq\pee(E_{n/2})\leq e^{-cn^{\frac{d-1}{d}}},\]
which is much smaller than $\pee(E_{u_n(x)})\pee(u_n(x)\leq |\cee (0)|<\infty)$.
Therefore, we obtain that for $y\neq z$
    \beq
    \label{FKG-applicsup}
    &&\pee(u_n(x)\leq |\ce(y)|<\infty, u_n(x)\leq |\ce (y)|<\infty)\\
    &&\qquad \leq
    \pee(E_{u_n(x)})\pee(u_n(x)\leq |\cee (0)|<\infty)(1+o(1)).\nonumber
    \eeq
We use (\ref{FKG-applicsup}), and thus
obtain
    \beq
    \pee\big(\t_{E_{u_n(x)}}\leq f_{E_{u_n(x)}}\big) &\geq& f_{E_{u_n(x)}}\pee(E_{u_n(x)})\\
    &&\qquad
    -f_{E_{u_n(x)}}^2\pee(E_{u_n(x)})\pee(u_n(x)\leq |\cee (0)|<\infty)(1+o(1)).
    \nonumber
    \eeq
Thus,
    \eq
    \lambda_{E_{u_n(x)}}\geq 1-f_{E_{u_n(x)}}
    \pee(u_n(x)\leq |\cee (0)|<\infty)(1+o(1))\geq 1-\pee(E_{u_n(x)})^{\beta}
    \en
for some $\beta>0$.
\epr

Finally, for supercritical clusters we expect
that
    \be\label{stret}
    \pee (n\leq |\ce(0)|<\infty) =C n^\alpha e^{-\eta n^{\delta}}[1+o(1)],
    \ee
i.e., Weibull tails (possibly with polynomial corrections),
and with $\delta=\frac{d-1}{d}$.

So far, (\ref{stret}) has not been proved rigorously, but
if we assume such a tail behavior, then we can
infer the precise form of the sequence $u_n(x)$ in
Lemma \ref{unx}.

\bp[Identification $u_n(x)$ for classical supercritical tails]
\label{suplem} If (\ref{stret}) is satisfied, then the sequence
$u_n(x)$ can be chosen of the form
\be
u_n(x) = \Big\lfloor \big(\frac{\log n}{\eta} + \frac{\alpha \log\log n}{\eta\delta} + x\big)^{1/\delta}
\Big\rfloor.
\ee
\ep

\bpr
Under the condition (\ref{stret}), it is a simple computation
to verify that
\be
\pee (u_n (x)\leq |\ce(0)|<\infty) = \frac{e^{-x}}{(2n+1)^d} (1+ o(1)).
\ee
\epr

\vskip0.5cm

\noindent
{\bf Proof of Theorem \ref{super} and Theorem \ref{thm-Gumbelsup}.}
We first finish the proof of Theorem \ref{super}. We follow the line of argument
in the proof of Theorem \ref{mainthm}.

We first use (\ref{Mntaueq}). Then we use Lemma \ref{lem-ttau} to
obtain that as long as $\pee (E_{u_n(x)})\leq n^{-d+\epsi}$, we have
(\ref{tothet}).

We again apply Theorem \ref{explaw}, and check the conditions. We
note from Lemma \ref{lem-locsub} that the events $E_{u_n(x)}$ are
localizable with local versions $E_{u_n(x)}'$. Furthermore, from
Proposition \ref{prop-secmomsup}, it follows that Condition
(\ref{secm}) is fulfilled for $E_{u_n(x)}$. Therefore, we may
apply Theorem \ref{explaw}.

We choose $u_n(x)$ as in Lemma \ref{unx}, and the event
$E_{u_n(x)}$ as before. Note that for this $u_n(x)$, we indeed
have that
    \eq
    \pee (E_{u_n(x)})=\frac{e^{-x}}{(2n+1)^d} a_n \leq n^{-d+\epsi},
    \en
so that we can use (\ref{tothet}).

Assume that $x\geq 0$. For $x<0$ some inequalities
reverse sign. Then we apply Theorem \ref{explaw} to obtain:
    \be
    \pee (\mee_n\geq u_n(x)) = \pee(\t_{E_{u_n(x)}}\leq (2n+1)^d)+o(1)=
    1-\exp\left(-\lambda_{E_{u_n(x)}}(2n+1)^d\pee (E_{u_n(x)})\right)+o(1).
    \ee
We need to investigate the exponent. By Lemma \ref{lem-taildistrsub}, we have that
    \eq
    \pee (E_{u_n(x)})
    = \frac{a_n}{(2n+1)^d}e^{-x}.
    \en
By Proposition \ref{prop-lambsup}, we have that
    \eq
    \lambda_{E_{u_n(x)}}=1+o(1).
    \en
Therefore, for any $x$,
    \be
    \pee (\mee_n\geq u_n(x))= 1- \exp(-a_n e^{-x}) + o(1).
    \label{conclsup}
    \ee
This completes the proof of Theorem \ref{super}.

If we further assume that $\pee$ has finite supercritical
clusters, then by Lemma \ref{unx} we can take $a_n=1+o(1)$. For Theorem
\ref{thm-Gumbelsup}, we note that the further assumption
(\ref{supass}) implies that $\pee$ has finite supercritical
clusters, and that with $u_n(x)=u_n+xu_n^{1/d}$, where
$u_n=u_n(0)$. Hence, we obtain that
    \eq
    \pee (E_{u_n(x)})=\pee (E_{u_n})\frac{\pee (E_{u_n+xu_n^{1/d}})}{\pee
    (E_{u_n})}
    =\pee (E_{u_n})e^{-x\eta\frac{d-1}{d}}[1+o(1)]
    =n^{-d}e^{-x\eta\frac{d-1}{d}}[1+o(1)].
    \en
The conclusion then follows from (\ref{conclsup}).
\QED

\subsection{Proof of Theorems \ref{cor-sub}, \ref{cor-sup} and \ref{thm-bc}}
\newcommand{\vep}{\varepsilon}
{\bf Proof of Theorems \ref{cor-sub} and \ref{cor-sup}.} We will
prove Theorems \ref{cor-sub} and \ref{cor-sup} simultaneously. In
order to do so, we let $\delta=1$ for $p<p_c$ and
$\delta=\frac{d-1}{d}$ for $p>p_c$. We then assume that
    \eq
    -\lim_{n\rightarrow \infty} \frac{1}{n^{\delta}} \log
    \pee(|\cee(0)|\geq n) = \xi
    \en
exists. The main ingredient is the following lemma:

    \bl[Convergence in probability]
    \label{lem-probbd}
    For any $\vep>0$, there exists $\kappa>0$ such that as $n\rightarrow \infty$,
        \eq
        \pee\big(\big|\frac{\mee_n}{(\log{n})^{1/\delta}}-C\big|>\vep)\leq
        n^{-\kappa},
        \en
    where $C=d\zeta$ for $p<p_c$ and $C=d^{\frac{d-1}{d}} \eta$
    for $p>p_c$.
    \el
     Before proving Lemma \ref{lem-probbd}, we will complete the
     proofs of Theorems \ref{cor-sub} and \ref{cor-sup} subject to
     Lemma \ref{lem-probbd}.

     Take $n_k=2^{k}$. As a consequence of  Lemma
     \ref{lem-probbd}, and the fact that for every $\kappa>0$,
        \[n_k^{-\kappa}=2^{-\kappa k}\]
     is summable in $k$, we obtain that
     $\frac{\mee_{n_k}}{(\log(n_k))^{1/\delta}}$
     converges to $C$ a.s. Thus, we have a.s.\ convergence along
     the subsequence $(n_k)_{k\geq 0}$. Moreover, we have that
     a.s.\ $n\mapsto \mee_n$ is non-decreasing. Therefore, for any $n_k<n\leq n_{k+1}$
     we can bound
        \eq\label{draak}
        \frac{\mee_{n_k}}{(\log(n_k))^{1/\delta}}\big(\frac{\log(n_k)}
        {\log(n_{k+1})}\big)^{1/\delta}\leq
        \frac{\mee_{n}}{(\log{n})^{1/\delta}}\leq
        \frac{\mee_{n_{k+1}}}{(\log(n_{k+1}))^{1/\delta}}
        \big(\frac{\log(n_{k+1})}{\log(n_k)}\big)^{1/\delta}.
        \en
    As $n\rightarrow \infty$, also $n_k, n_{k+1}\rightarrow
    \infty$. Thus, $\frac{\mee_{n_{k}}}{(\log(n_{k}))^{1/\delta}}$ and
    $\frac{\mee_{n_{k+1}}}{(\log(n_{k+1}))^{1/\delta}}$ converge a.s.\ to
    $C$. Furthermore,
        \eq
        \lim_{k\rightarrow \infty} \frac{\log(n_{k+1})}
        {\log(n_{k})}=\lim_{k\rightarrow \infty} \frac{k+1}{k}=1,
        \en
so that both upper and lower bound in (\ref{draak}) converge to $C$ almost surely.
This completes the proofs of Theorems \ref{cor-sub} and \ref{cor-sup}.
    \qed
    \vskip0.5cm

    \noindent
    {\it Proof of Lemma \ref{lem-probbd}.}  Fix $\vep >0$. We will prove
        \eq
        \pee\Big(\frac{\mee_n}{(\log{n})^{1/\delta}}>C+\vep\Big)
        \leq n^{-\kappa},
        \label{bdkappa1}
        \en
    and
        \eq
        \pee\Big(\frac{\mee_n}{(\log{n})^{1/\delta}}<C-\vep\Big)
        \leq n^{-\kappa}.
        \label{bdkappa2}
        \en
    Let $C$ be the constant
    such that
        \eq
        \label{krawaat}
        \pee\Big(\frac{|\cee(0)|}{(\log{n})^{1/\delta}}>C\Big)
        =n^{-d(1+o(1))},
        \en
This constant exists by (\ref{zetadef}) in the case $p<p_c$, and
by assumption (\ref{etadef}) (proved in $d=2,3$) for $p>p_c$.

    With this choice of $C$, for $\vep>0$, there exists a $\kappa' \in (0, d)$ such that
        \eq
        \pee\Big(\frac{|\cee(0)|}{(\log{n})^{1/\delta}}>C+\vep\Big)
        \leq n^{-d-\kappa'},
        \label{bd1}
        \en
    while
        \eq
        \pee\Big(\frac{|\cee(0)|}{(\log{n})^{1/\delta}}<C-\vep\Big)
        \leq 1-n^{-d+\kappa'}.
        \label{bd2}
        \en

    To prove (\ref{bdkappa1}), we use that
        \beq
        \pee\Big(\frac{\mee_n}{(\log{n})^{1/\delta}}>C+\vep\Big)
        &=&\pee\Big(\bigcup_{x\in B_n}
        \{|\cee(x)|>(C+\vep)(\log{n})^{1/\delta}\}\Big)\nonumber\\
        &\leq& \sum_{x\in B_n}
        \pee\Big(|\cee(0)|>(C+\vep)(\log{n})^{1/\delta}\Big)\nonumber\\
        &\leq& |B_n| n^{-d-\kappa'}\leq n^{-\kappa},
        \label{bdkappa1a}
        \eeq
    where we use (\ref{bd1}).

    To prove (\ref{bdkappa2}), we use that the events
    $\{|\cee(x)|\leq (C+\vep)(\log{n})^{1/\delta}\}_{x\in A_n}$ are
    independent when
    \eq
    A_n= (K_n\Z)^d \cap B_n.
    \en
    and
    \be
    K_n= \lceil
    (C+\vep)(\log{n})^{1/\delta}\rceil
    \en
    Therefore,
        \beq
        \pee\Big(\frac{\mee_n}{(\log{n})^{1/\delta}}<C-\vep\Big)
        &=&\pee\Big(\bigcap_{x\in A_n} \{|\cee(x)|\leq (C-\vep)(\log{n})^{1/\delta}\}\Big)\nonumber\\
        &=& \prod_{x\in A_n} \pee\Big(|\cee(x)|<(C-\vep)(\log{n})^{1/\delta}\Big)\nonumber\\
        &\leq&
        \pee\Big(|\cee(0)|<(C-\vep)(\log{n})^{1/\delta}\Big)^{|A_n|}.
        \label{bdkappa2a}
        \eeq
    We next use (\ref{bd2}) and the fact that
        \eq
        |A_n| \geq \Big(\frac{n}{\lceil
        (C+\vep)(\log{n})^{1/\delta}\rceil}\Big)^d,
        \en
    so arrive at a bound, for every $\kappa\in (0,\kappa')$,
         \beq
        \pee\Big(\frac{\mee_n}{(\log{n})^{1/\delta}}<C-\vep\Big)
        &\leq&
        \big(1-n^{-d+\kappa'}\big)^{|A_n|}\leq n^{-\kappa},
        \label{bdkappa2b}
        \eeq
    which completes the proof.
    \qed

    \vskip 0.5cm

    \noindent
    {\bf Proof of Theorem \ref{thm-bc}.} We again use
    (\ref{krawaat}) together with the observation that
    the events $\{\mee_n^{\scriptscriptstyle {\rm (zb)}}
    \neq \mee_n^{\scriptscriptstyle {\rm (fb)}}\}$ and
    $\{\mee_n^{\scriptscriptstyle {\rm (zb)}}
    \neq \mee_n^{\scriptscriptstyle {\rm (pb)}}\}$
    are contained in the event that there exists a cluster
    on the boundary (either with free or periodic boundary
    conditions) such that there exists an $x\in \partial B_n$
    such that $|\cee(x)|\geq \mee_n^{\scriptscriptstyle {\rm (zb)}}$.
    By Theorems \ref{cor-sub} and \ref{cor-sup}, we have that
    $\mee_n^{\scriptscriptstyle {\rm (zb)}}\geq
    (C-\vep)(\log{n})^{1/\delta}$ a.s.
    By (\ref{krawaat}) and when $\vep>0$ is sufficiently small,
    this probability is thus bounded above by
        \[
        n^{d-1}\pee\Big(|\cee(x)|\geq
        (C-\vep)(\log{n})^{1/\delta}\Big) \leq n^{-\kappa}
        \]
    for some $\kappa>0$.
    \qed

\section*{Acknowledgements}
The work of RvdH was supported in part by Netherlands Organisation for
Scientific Research (NWO).
The authors thank Harry Kesten for stimulating discussions during
the early stages of this work.

\end{document}